\documentclass{amsart}
\usepackage{amssymb,latexsym}
\usepackage{amscd}
\usepackage[all]{xy}

\newtheorem{theorem}{Theorem}[section]

\newtheorem{prop}[theorem]{Proposition}

\newtheorem{coro}[theorem]{Corollary}

\newcommand{\ext}{\mbox{\rm Ext}}
\theoremstyle{remark}
\newtheorem{remark}[theorem]{Remark}
\renewcommand{\dim}{\mbox{\rm dim}}

\newcommand{\im}{\mbox{\rm Im}}

\renewcommand{\hom}{\mbox{\rm Hom}}

\newcommand{\cor}{\mbox{\rm cor}}
\newcommand{\res}{\mbox{\rm res}}
\newcommand{\bZ}{\mathbb{Z}}
\newcommand{\bQ}{\mathbb{Q}}
\newcommand{\sK}{\mathcal{K}}

\newcommand{\fS}{\mathfrak{S}}
\newcommand{\bF}{\mathbb{F}}
\newcommand{\bR}{\mathbb{R}}
\newcommand{\bC}{\mathbb{C}}
\newcommand{\gal}{\mbox{\rm Gal}}
\newcommand{\rk}{\mbox{\rm Rank}}
\newcommand{\len}{\mbox{\rm length}}
\numberwithin{equation}{section}

\newcommand{\abs}[1]{\lvert#1\rvert}

\begin{document}

\title{Introduction to Local and Global Euler Characteristic Formulas}

\author{Wei Lu}
\address{Department of Mathematical Science, Tsinghua University, Beijing, P. R. China, 100084}

\email{weillu19880909@gmail.com}

%    General info
\date{November 30, 2011 and, in revised form, December 7, 2011.}

\begin{abstract}
This is a note of talks I gave at the number theory seminar at Tsinghua University in Fall 2011.

We will introduce the local and global Euler characteristic formulas given by John Tate(1962) for Galois cohomology.
We will give a detailed proof based on the idea in Hida's book[1,Ch4.4.4 and 4.4.5] and Milne's book[2,Ch1.5].
\end{abstract}

\maketitle
This note is organized as follows. In preliminary, we review the definition of group cohomology and some basic properties. In Section 1, we will give a detailed proof of the local case. In Section 2, we also prove the global case by using a powerful theorem given by John Tate. Both of these proofs roughly follow Hida's book[1, Ch4.4.4 and 4.4.5] and Milne's book[2,Ch1.5].
\section*{Preliminary}
In the section, we recall some basic facts on Galois cohomology without proof. Readers can see [4] or [5] for more details.

\subsection*{Cohomology and Cochains}

Let $G$ be a finite group, and $M$ be a $G$-module. The functor $M\mapsto M^G$ from the category of $G$-modules to the category of Abelian groups is left-exact. The derived functor is denoted by $H^n(G,-)$. $H^n(G,M)$ is called the cohomology group, and it can be computed by the complex of cochains in the following way.

We define $C^n(G,M):=Map(G^n,M)$, an element of $C^n(G,M)$ is a function $f$ of $n$ variables in $G$ with codomain $M$.
The differential maps
$$d_n:C^{n}(G,M)\rightarrow C^{n+1}(G,M)$$
are defined by
\begin{eqnarray}
d_n(f_n)(g_1,\cdots,g_n)&=&g_1\cdot f(g_2,\cdots,g_n)+\sum_{i=1}^n(-1)^if(g_1,\cdots,g_ig_{i+1},\cdots,g_{n+1})\nonumber\\
&&+(-1)^{n+1}f(g_1,\cdots,g_{n}),\nonumber
\end{eqnarray}
and we can check $d_{n+1}\circ d_n=0$ directly. The cohomological groups are given by $$H^n(G,M)=H^n(C^\bullet(G,M))=\ker d_n/\im d_{n-1}.$$

When $G$ is a profinite group, the discrete abelian groups on which $G$ acts continuously form an abelian category $C_G$, which is a full subcategory of the category of all $G$-modules. For a (discrete) $G$-module $M$, we define
$$H^n(G,M):=\varinjlim\limits_{H\triangleleft G}H^n(G/H,M^H),$$
where $H$ runs over the open normal subgroups of $G$. If we want to compute the cohomology groups $H^n(G,M)$, we can still use the same method applied to compute the cochains. The only change is that the cochains must be continuous.

\subsection*{cor$_{G/U}$, res$_{G/U}$}
Let $U$ be a closed subgroup of $G$, we have restriction map $$\res_{G/U}:H^n(G,M)\rightarrow H^n(U,M).$$ If $U$ is a open subgroup of $G$ of finite index, we have corestriction map $$\cor_{G/U}:H^n(U,M)\rightarrow H^n(G,M).$$

Then we have the following proposition.
\begin{prop}The following properties hold:
\begin{enumerate}
  \item $\cor_{G/U}\circ\res_{G/U}(x)=[G:U]x$, where $x\in H^n(G,M)$;
  \item If $G$ is a finite group, $M$ is a finite $G$-module, and $(\abs{G},\abs{M})=1$, then $H^q(G,M)=0$ for all $q>0$.
\end{enumerate}
\end{prop}

\subsection*{Inflation and Restriction Sequences}

\begin{prop}
Let $U$ be a closed normal subgroup of $G$, and suppose $H^p(U,M)=0$ for all $p=1,2,\cdots,q-1$. Then the following sequence is exact:
$$0\rightarrow H^q(G/U,M^U)\rightarrow H^q(G,M)\rightarrow H^0(G/U,H^q(U,M))\rightarrow H^{q+1}(G/U,M^U).$$
\end{prop}
\subsection*{Shapiro's Lemma} Let $U$ be a closed subgroup of $G$. The induced module is given by $$Ind_U^GM=\hom_{\bZ[U]}(\bZ[G],M).$$
\begin{prop}
We have isomorphisms:
$$H^n(G,Ind_U^GM)=H^n(U,M)$$
for all $n\geq0$.

\end{prop}
\subsection*{Tate Cohomology}
Let $G$ be a finite group and $M$ be a $G$-module, then the Tate groups are defined by:
$$H^n_T(G,M)=\left\{
               \begin{array}{ll}
                 H^n(G,M), & n\geq1 \\
                 H^0(G,M), & n=0 \\
                 \ker N_m/I_G M, & n=-1 \\
                 H_{i-1}(G,M), & n<-1
               \end{array}
             \right.
$$
where $N_m:m\mapsto\sum\limits_{g\in G}gm$, $I_G\subset\bZ[G]$ is generated by $-1$.
\begin{prop}
Suppose that $G$ is a finite cyclic group generated by $g$. Then
$$H^{n}_T(G,M)\cong H^{n+2}_T(G,M).$$
\end{prop}
\begin{coro}
All the notations are the same as above, then we have:
$$H^{2n}(G,M)\cong H^0_T(G,M)=M^G/N_m(M)$$
for all $n\geq1$, and
$$H^{2n-1}(G,M)\cong\ker N_m/I_GM.$$
\end{coro}
\begin{prop}
If $G$ is cyclic, and $M$ is finite, then
$$|H^0_T(G,M)|=|H^1_T(G,M)|.$$
\end{prop}

\section{Local Euler Characteristic Formula}
The main result is the following theorem. \begin{theorem}[Local Case]
Let $K/\mathbb{Q}_p$ be a finite extension for a prime $p$, $G=Gal(\overline{\mathbb{Q}_p}/K)$ and M be a finite (discrete) $G$-module. We have local Euler characteristic formula:
$$\frac{|H^0(G,M)|\cdot|H^2(G,M)|}{|H^1(G,M)|}=\frac{|H^0(G,M)|\cdot|H^0(G,M^*(1))|}{|H^1(G,M)|}=|| M||_K,$$
where $M^*(1)=\hom(M,\overline{K}^\times)$, $|n|_K=[O_K:nO_K]^{-1}$ for a positive integer $n$.
\end{theorem}

The first equation is as a result of Tate duality.
\begin{prop}
(Tate duality) Let $M$ be a finitely generated discrete $\bZ[G]$-module,
$$H^r(G,M^*(1))\cong H^{2-r}(G,M)^*  $$
for all $0\leq r\leq2$, where $M^*(1)=\hom(M,\overline{K}^\times)$, $N^*=\hom(N,\mathbb{Q}/\mathbb{Z})$ is the Potryagin dual module of an abelian group $N$.

In particular, if $M$ is finite, all cohomology groups $H^r(G,M)$ are finite and $H^r(G,M)=0$ for $r\geq3$.
\end{prop}

\subsection*{Proof of the Local Case}
We simply write $H^n(M)$ for $H^n(G,M)$. Since $M=\bigoplus\limits_{l}M[l^\infty]$ for prime $l$, so we only need to prove the case $M=M[l^\infty]$ because of additions. Now we may assume that $M=M[l^\infty]$, then $H^q(M)$ is a $\bZ_l$-module of finite length. For any finite $\bZ_l$-module $N$, we have $\abs{N}=l^{\len_{\bZ_l}(N)}$ because a simple non-zero $\bZ_l$-module must be isomorphic to $\bZ_l/l$.  Here the $\len(N)$ is the length of the Jordan-Holder sequence of $\bZ_l$-module $M$.

We define the local Euler character by
$$\chi(M)=\chi(G,M)=\sum\limits_{q=0}^2(-1)^q\len_{\bZ_l}H^q(M),$$
$$\chi'(M)=\chi'(G,M)=\log_l(\abs{|M|}_K)=\left\{
                                        \begin{array}{ll}
                                          -[K:\bQ_p]\len_{\bZ_p}M & \hbox{$l=p$}, \\
                                          0 & \hbox{$l\neq p$}.
                                        \end{array}
                                      \right.
$$
Note that the left and the right side of this formula, so we only need to prove $$\chi(M)=\chi'(M)=\left\{
                                        \begin{array}{ll}
                                          -[K:\bQ_p]\len_{\bZ_p}M & \hbox{$l=p$}, \\
                                          0 & \hbox{$l\neq p$}.
                                        \end{array}
                                      \right.$$

We first check the formula for the trivial case: $M=\bF_l=\bZ/l\bZ$ ($G$ acts on $\bF_l$ trivially). By Tate Duality, $\dim_{\bF_l}H^0(G,\bF_l)=\dim_{\bF_l}\bF_l^G=1$, $\dim_{\bF_l}H^2(G,\bF_l)=\dim_{\bF_l}H^0(G,\mu_l)^*=\dim_{\bF_l}\mu_l^*(K)=\dim_{\bF_l}\mu_l(K)$, where $\mu_l(K)=\{z\in K|z^l=1\}$. On the other hand, by Kummer theory, $\dim_{\bF_l}H^1(G,\bF_l)=\dim_{\bF_l}H^1(G,\mu_l)^*=\dim_{\bF_l}H^1(G,\mu_l)=\dim_{\bF_l}K^\times/(K^\times)^l$. The reason for the last equation is that
\[\begin{CD}
1@>>>\mu_l@>>>\overline{K}^\times@>x\mapsto x^l>>\overline{K}^\times@>>>1.
\end{CD}\]
Then we have
\[\begin{CD}
1@>>>H^0(G,\mu_l)@>>>H^0(G,\overline{K}^\times)@>>>H^0(G,\overline{K}^\times)@>>>H^1(G,\mu_l)@>>>H^1(G,\overline{K}^\times)\\
@.@|@|@|@|@|\\
1@>>>\mu_l(K)@>>>K^\times@>x\mapsto x^l>>K^\times@>>>H^1(G,\mu_l)@>>>1,
\end{CD}\]
where $H^1(G,\overline{K}^\times)=1$ is given by Hilbert 90.
So $H^1(G,\mu_l)\cong K^\times/(K^\times)^l$. Since $K^\times\cong O_K^\times\times\bZ$, $O^\times\cong O_K\times\mu\cong O_K\times\mu_{l^\infty}(K)\times\prod\limits_{q\neq l}\mu_{q^\infty}(K)$, where $O_K$ is the integer ring of $K$. So
$$K^\times/(K^\times)^l\cong\left\{
                              \begin{array}{ll}
                                \bZ/l\bZ\oplus\mu_l(K) & l\neq p, \\
                                \bZ/p\bZ\oplus O_K/pO_K\oplus\mu_p(K) & l=p.
                              \end{array}
                            \right.
$$
When $l=p$,
\begin{eqnarray}
\chi(\bF_p)&=&\dim_{\bF_p}H^0(G,\bF_p)-\dim_{\bF_p}H^1(G,\bF_p)+\dim_{\bF_p}H^2(G,\bF_p)\nonumber\\
&=&1-\dim_{\bF_p}(\bF_p\oplus O_K/pO_K\oplus\mu_p(K))+\dim_{\bF_p}\mu_p(K)\nonumber\\
&=&-\dim_{\bF_p}O_K/pO_K=-[K:\bQ_p]\dim_{\bF_p}\bF_p=\chi'(\bF_p).\nonumber
\end{eqnarray}
And when
$l\neq p$, $$\chi(\bF_l)=1-\dim_{\bF_l}(\bF_l\oplus\mu_l(K))-\dim_{\bF_l}\mu_l(K)=0=\chi'(\bF_l).$$
By Tate duality, the formula also holds for $M=\mu_l=\bF_l^*(1)$.

We will prove that if $0\rightarrow L\rightarrow M\rightarrow N\rightarrow 0$ is an exact sequence of finite $\bZ_l[G]$-modules, then$\chi'(M)=\chi'(N)+\chi'(L)$, $\chi(M)=\chi(N)+\chi(L)$, and $\chi'(M)=\chi'(M^{ss})$, $\chi(M)=\chi(M^{ss})$, where $M^{ss}=\bigoplus\limits_{q=1}^nM_q/M_{q-1}$ for a Jordan-Holder sequence $0=M_0\subset M_1\subset \cdots M_n=M$ of $\bZ_l[G]$-modules.

In fact, $\chi'(M)=-[K:\bQ_p]\len_{\bZ_p}M=-[K:\bQ_p](\len_{\bZ_p}L+\len_{\bZ_p}N)=\chi'(N)+\chi'(L)$. And $0\rightarrow H^0(L)\rightarrow H^0(M)\rightarrow H^0(N)\rightarrow H^1(L)\rightarrow H^1(M)\rightarrow H^1(N)\rightarrow H^2(L)\rightarrow H^2(M)\rightarrow  H^2(N)\rightarrow 0$ is an exact sequence which follows Proposition 1.2. So we have
\begin{eqnarray}
\chi(M)&=&\sum\limits_{q=0}^2(-1)^q\len H^q(M)\nonumber\\
&=&\sum\limits_{q=0}^2(-1)^q(\len H^q(L)+\len H^q(N))\nonumber\\
&=&\sum\limits_{q=0}^2(-1)^q\len H^q(L)+\sum\limits_{q=0}^2(-1)^q\len H^q(N)\nonumber\\
&=&\chi(L)+\chi(N).\nonumber
\end{eqnarray}
Thus we get
\begin{eqnarray}
\chi'(M^{ss})&=&\chi'(\bigoplus\limits_{q=1}^nM_q/M_{q-1})\nonumber\\
&=&\sum\limits_{q=1}^n\chi'(M_q/M_{q-1})\nonumber\\
&=&\sum\limits_{q=1}^n(\chi'(M_q)-\chi'(M_{q-1}))\nonumber\\
&=&\chi'(M_n)=\chi'(M).\nonumber
\end{eqnarray}
And the same reason for $\chi(M)=\chi(M^{ss})$.

However, $M^{ss}$ is a $\bF_l[G]$-module because $M_q/M_{q-1}$ is a $\bF_l[G]$-module. Then we may assume that $M$ itself is a $\bF_l[G]$-module. At this time $\dim_{\bF_l}M=\len_{\bZ_l}M$.

Now we recall the notation of Grothendieck groups. Let $G$ be a profinite group and $E$ be a field. We consider the category $Rep_E(G)$ made up of the following data:
\begin{enumerate}
  \item Objects are finite dimensional $E-$vector spaces with a continuous action of $G$ under the discrete topology;
  \item Morphisms are $E[G]-$linear maps.
\end{enumerate}
Grothendieck group $R_E(G)$ of $Rep_E(G)$ is an Abelian group which is defined by generators and relations: $R_E(G)$ is generated by symbols $[M]$ for objects $M\in Rep_E(G)$. The only relation is $[M]=[N]+[L]$ if $0\rightarrow L\rightarrow M\rightarrow N\rightarrow 0$ is a short exact sequence.

Now we consider the category $Rep_{\bF_l}(G)$ which is made of all finite $\bF_l[G]$-module. Its Grothendieck group is $R_{\bF_l}(G)$, we can regard $\chi$ and $\chi'$ as functions on the Grothendieck group $R_{\bF_l}(G)$ with value $\bZ$. We need to check the formula for a set of generators of $R_{\bF_l}(G)$. As $\bZ$ is torsion-free, we only check it for a set of generators of $R_{\bF_l}(G)\otimes_{\bZ}\bQ$. And we can find a set of generators by the following proposition:

\begin{prop}
(see $[2,\mathrm{lemma} 2.10]$) Let $G$ be a finite group and, for any subgroup H of G, let $Ind^G_H$ be the homomorphism $R_{\bF_p}(H)\otimes\bQ\rightarrow R_{\bF_p}(G)\otimes\bQ$ taking the class of an $H$-module
to the class of the corresponding induced $G$-module. $R_{\bF_p}(G)\otimes\bQ$ is generated by the images of the $Ind^G_H$ as $H$ runs over the set of cyclic subgroups of $G$ of order prime to $p$.
\end{prop}

We take a finite Galois extension $F/K$ such that $\gal(\overline{\bQ}_p/F)$ acts trivially on $M$, write $\overline{G}=\gal(F/K)$. Hence, We only need to check the formula for a set of generators of $R_{\bF_l}(\overline{G})$. However, by Proposition 1.4, $R_{\bF_l}(\overline{G})\otimes_{\bZ}\bQ$ is generated by $Ind_{\overline{H}}^{\overline{G}}\rho$ for cyclic subgroups $\overline{H}$ of order prime to $l$ and character $\rho:\overline{H}\rightarrow \sK^\times$ for a finite extension $\sK/\bF_l$. Thus we can assume $M=Ind_{\overline{H}}^{\overline{G}}\rho=Ind_{H}^{G}\rho$, where $H=\gal(\overline{\bQ}_p/F^{\overline{H}}).$ Then by Shapiro's lemma, $H^q(G,Ind_{\overline{H}}^{\overline{G}}\rho)\cong H^q(G,Ind_{H}^{G}\rho)\cong H^q(H,\rho)$, thus $\chi(G,Ind_{H}^{G}\rho)=\chi(H,\rho)$. So we only need to check the formula for $\rho$ (or for one-dimensional single module $V(\rho)$ on which $\overline{H}$ acts via $\rho$). \\
$$
%\xymatrix{\overline{\bQ}_p  \ar@{-}[d] \ar@/_/[d]_{G'}\ar@/^10mm/[ddd]^{G}\\
%F \ar@{-}[d] \ar@/_/[dd]_{\overline{G}} \ar@/^/[d]^{\overline{H}}\\
%F^{\overline{H}} \ar@{-}[d]\\
%K}
$$
Thus we may assume $F^{\overline{H}}=K$, then $\overline{H}=\overline{G}$, $M=V(\rho)$ is one-dimensional over $\sK$ and $\overline{G}=\overline{H}$ is a cyclic group, $(|\overline{G}|,l)=1$. By Proposition 0.1(2), we know $H^q(\overline{G},M)=0$ for all $q\geq 1$. Hence, we have inflation and restriction sequence,
\[\begin{CD}
H^q(G/G',M^{G'})@>>>H^q(G,M)@>>>H^0(G/G',H^q(G',M))@>>>H^{q+1}(G/G',M^{G'})\\
@|@|@|@| \\
0=H^q(\overline{G},M)@>>>H^q(G,M)@>>>H^0(\overline{G},H^q(G',M))@>>>H^{q+1}(\overline{G},M)=0,
\end{CD}\]
where $G'=\gal(\overline{Q}_p/F)$. Moreover, we have $H^q(G,M)\cong H^0(\overline{G},H^q(G'M))$ for $q=0,1,2$, and we note that
$$H^q(G',M)=H^q(G',\sK)=H^q(G',\bF_l)\otimes_{\bF_l}\sK=\left\{
                                                       \begin{array}{ll}
                                                         \sK & q=0;\\
                                                         ((F^\times/(F^\times)^l )^*\otimes_{\bF_l}\sK& q=1;\\
                                                         \mu_l^*(F)\otimes_{\bF_l}\sK & q=2.
                                                       \end{array}
                                                     \right.
$$
Then $$\chi(G,M)=\dim_{\bF_l}\sK^{\overline{G}}-\dim_{\bF_l}(((F^\times)/(F^\times)^l)^*\otimes_{\bF_l}\sK)^{\overline{G}}+\dim_{\bF_l}(\mu_l^*(F)\otimes_{\bF_l}\sK)^{\overline{G}}.$$ Since we have checked the cases $M=\bF_l$ and $M=\mu_l$, we may assume that $\rho$ is neither trivial nor cyclotomic character . Hence $\sK^{\overline{G}}=(\mu_l^*(F)\otimes_{\bF_l}\sK)^{\overline{G}}=0$ because the action of the Galois group on $\sK$ is via nontrivial character $\rho$ on $\mu_p$ is via the cyclotomic character. Therefore, $\chi(G,M)=-\dim_{\bF_l}(((F^\times)/(F^\times)^l)^*\otimes_{\bF_l}\sK)^{\overline{G}}$.

When $l=p$, we need to show that
$$\chi(G,M)=-\dim_{\bF_l}(((F^\times)/(F^\times)^l)^*\otimes_{\bF_l}\sK)^{\overline{G}}=-[K:\bQ_p]\dim_{\bF_p}M=\chi'(G,M).$$
Since $\mu$ is the maximal torsion-subgroup of $F^\times$, then we have
\[\begin{CD}
1@>>>\mu@>>>F^\times@>>>F^\times/\mu@>>>1\\
@.@VVpV@VVpV@VVpV\\
1@>>>\mu@>>>F^\times@>>>F^\times/\mu@>>>1,
\end{CD}\]
where $p:x\mapsto x^p$. By snake Lemma, we have
$$1\rightarrow\mu/\mu^p\rightarrow (F^\times)/(F^\times)^p\rightarrow(F^\times/\mu)/(F^\times/\mu)^p\rightarrow1.$$
As "$^*$" is defined as following: "$^*"=\hom(-,\bQ/\bZ)=\hom(-,\bF_p)$ (only in this situation) is a contravariant and left exact functor, thus
$$1\rightarrow((F^\times/\mu)/(F^\times/\mu)^p)^*\rightarrow ((F^\times)/(F^\times)^p)^*\rightarrow(\mu/\mu^p)^*\rightarrow\ext^1((F^\times/\mu)/(F^\times/\mu)^p,\bF_p)=1.$$

We know that $\sK$ is flat, then
$$1\rightarrow((F^\times/\mu)/(F^\times/\mu)^p)^*\otimes_{\bF_p}\sK\rightarrow ((F^\times)/(F^\times)^p)^*\otimes_{\bF_p}\sK\rightarrow(\mu/\mu^p)^*\otimes_{\bF_p}\sK\rightarrow1.$$

Hence
$$1\rightarrow H^0(\overline{G},((F^\times/\mu)/(F^\times/\mu)^p)^*\otimes_{\bF_p}\sK)\rightarrow H^0(\overline{G},((F^\times)/(F^\times)^p)^*\otimes_{\bF_p}\sK)\rightarrow H^0(\overline{G},(\mu/\mu^p)^*\otimes_{\bF_p}\sK).$$

However, $\mu/\mu^p\cong\mu_p(F)$, then $$H^0(\overline{G},(\mu/\mu^p)^*\otimes_{\bF_p}\sK)=(\mu_p^*(F)\otimes_{\bF_p}\sK)^{\overline{G}}=0.$$

Therefore,
$$H^0(\overline{G},((F^\times/\mu)/(F^\times/\mu)^p)^*\otimes_{\bF_p}\sK)\cong H^0(\overline{G},((F^\times)/(F^\times)^p)^*\otimes_{\bF_p}\sK).$$

In other words,
$$\dim_{\bF_p}(((F^\times)/(F^\times)^p)^*\otimes_{\bF_p}\sK)^{\overline{G}}=\dim_{\bF_p}(((F^\times/\mu)/(F^\times/\mu)^p)^*\otimes_{\bF_p}\sK)^{\overline{G}}=\dim_{\bF_p}(((F^\times/\mu)\otimes_{\bZ}\bF_p)^*\otimes_{\bF_p}\sK)^{\overline{G}}.$$

Writing the additive valuation of $F$ as $v:\bF^\times\rightarrow\bZ$, then we have an exact sequence:
$$1\rightarrow O_F^\times/\mu\rightarrow F^\times/\mu\rightarrow\bZ\rightarrow0.$$
Then the exact sequence is torsion-free, and after tensor $\bF_p$, we still have an exact sequence
$$1\rightarrow O_F^\times/\mu\otimes_{\bZ}\bF_p\rightarrow F^\times/\mu\otimes_{\bZ}\bF_p\rightarrow\bF_p\rightarrow0.$$

Using $*$ functor again, we have:
$$0\rightarrow\hom(\bF_p,\bF_p)=\bF_p\rightarrow(\bF^\times/\mu\otimes_{\bZ}\bF_p)^*\rightarrow(O_F^\times/\mu\otimes_{\bZ}\bF_p)^*\rightarrow\ext^1(\bF_p,\bF_p)=0.$$
So we have $$0\rightarrow\sK^{\overline{G}}\rightarrow((\bF^\times/\mu\otimes_{\bZ}\bF_p)^*\otimes_{\bF_p}\sK)^{\overline{G}}\rightarrow((O_F^\times/\mu\otimes_{\bZ}\bF_p)^*\otimes_{\bF_p}\sK)^{\overline{G}}\rightarrow H^1(\overline{G},\sK)=0.$$
Therefore,
$$
\dim_{\bF_p}((F^\times/\mu\otimes_{\bZ}\bF_p)^*\otimes_{\bF_p}\sK)^{\overline{G}}=\dim_{\bF_p}((O_{F}^\times/\mu\otimes_{\bZ}\bF_p)^*\otimes_{\bF_p}\sK)^{\overline{G}}.
$$

Now we want to lift the representation $\rho$ to characteristic $0$ representation $\widetilde{\rho}$ by the following proposition.
\begin{prop}
(see [1, corollary 2.7]) Let $K$ be a finite extension of $\overline{\bQ}_p$ with $p-$adic integer ring $O$. Let $E=O/m_O$ for the maximal ideal $m_O$ of $O$. Suppose $p| |G|$ is not true and that all irreducible representations of $G$ over $K$ are absolutely irreducible. Then all irreducible representations of $G$ over $E$ are absolutely irreducible, and the reduction map $\rho\mapsto(\rho \mod m_O)$ induces a bijection between isomorphism classes of absolutely irreducible representations of $G$ over $K$ and over $E$, preserving dimension.
\end{prop}

 For that we take the unique unramified extension $L$ of $\bQ_p$ of degree $\dim_{\bF_p}\sK$. Then we have $O_L/(p)\cong\sK$, $O^\times_L\cong(1+pO_L)\times\sK^\times$, where $O_L$ is the p-adic integer ring of $L$. By the isomorphism, we may think $\rho$ has valuation in $O_L^\times$. We write the character $\widetilde{\rho}:\overline{G}\rightarrow O_L^\times$, which is called the Teichmuller lift of $\rho$. Since $O_F^\times/\mu$ is torsion free and $(|\overline{G}|,\rho)=1$, by Proposition 1.4 for the unique Teichmuller lift $\widetilde{\rho}$ of $\rho$, we have:
\begin{eqnarray}
\dim_{\bF_p}((O_{F}^\times/\mu\otimes_{\bZ}\bF_p)^*\otimes_{\bF_p}\sK)^{\overline{G}}&=&\dim_{\bF_p}((O_{F}^\times/\mu\otimes_{\bZ_p}\bF_p)^*\otimes_{\bF_p}\sK)^{\overline{G}}\nonumber\\
&=&\dim_{\bF_p}(\hom(O_{F}^\times/\mu\otimes_{\bZ_p}\bF_p,\bF_p)\otimes_{\bF_p}\sK)^{\overline{G}}\nonumber\\
&=&\rk_{\bZ_p}(\hom(O_{F}^\times/\mu\otimes_{\bZ_p}\bZ_p,\bZ_p)\otimes_{\bZ_p}O_L)^{\overline{G}}\nonumber\\
&=&\dim_{\bQ_p}(\hom(O_{F}^\times/\mu\otimes_{\bZ}\bQ,\bQ_p)\otimes_{\bQ_p}L)^{\overline{G}}.\nonumber
\end{eqnarray}
By p-adic logarithm, we know that $O_F^\times/\mu\otimes_{\bZ}\bQ\cong F$ as $\overline{G}$-module. Hence
$$\dim_{\bQ_p}(\hom((O_{F}^\times/\mu)\otimes_{\bZ}\bQ,\bQ_p)\otimes_{\bQ_p}L)^{\overline{G}}=\dim_{\bQ_p}(\hom(\bF,\bQ_p)\otimes_{\bQ_p}L)^{\overline{G}}.$$
By normal base theorem, $F\cong K[\overline{G}]\cong \bQ_p[\overline{G}]^{[K:\bQ_p]}$. Then
\begin{eqnarray}
\dim_{\bQ_p}(\hom(F,\bQ_p)\otimes_{\bQ_p}L)^{\overline{G}}&=&\dim_{\bQ_p}(\hom(\bQ_p[\overline{G}]^{[K:\bQ_p]},\bQ_p)\otimes_{\bQ_p}L)^{\overline{G}}\nonumber\\
&=&[K:\bQ_p]\dim_{\bQ_p}(\hom(\bQ_p[\overline{G}],\bQ_p)\otimes_{\bQ_p}L)^{\overline{G}}.\nonumber
\end{eqnarray}

We can easily check that $\hom(\bQ_p[\overline{G}],\bQ_p)\cong\bQ_p[\overline{G}]$ as $\overline{G}$-module by $\psi:f\mapsto \sum\limits_{\sigma\in \overline{G}}a_{\sigma}\sigma$, where $a_{\sigma}=f(\sigma)$, then
\begin{eqnarray*}
& &[K:\bQ_p]\dim_{\bQ_p}(\hom(\bQ_p[\overline{G}],\bQ_p)\otimes_{\bQ_p}L)^{\overline{G}}\nonumber\\
&=&[K:\bQ_p]\dim_{\bQ_p}(\bQ_p[\overline{G}]\otimes_{\bQ_p}L)^{\overline{G}}\nonumber\\
&=&[K:\bQ_p]\dim_{\bQ_p}(\bQ_p[\overline{G}]\otimes_{\bQ_p}L^0)^{\overline{G}}\cdots\cdots(*)\\
&=&[K:\bQ_p]\dim_{\bQ_p}(\bQ_p[\overline{G}]^{\overline{G}}\otimes_{\bQ_p}L^0)\nonumber\\
&=&[K:\bQ_p]\dim_{\bQ_p}(\bQ_p\otimes_{\bQ_p}L^0)\nonumber\\
&=&[K:\bQ_p]\dim_{\bQ_p}L^0\nonumber\\
&=&[K:\bQ_p]\dim_{\bQ_p}L\nonumber\\
&=&[K:\bQ_p]\dim_{\bF_p}\sK\nonumber\\
&=&[K:\bQ_p]\dim_{\bF_p}M,\nonumber
\end{eqnarray*}
where $\bQ_p[\overline{G}]^{\overline{G}}\cong\bQ_p$ and (*) follows from the isomorphism:
$$\bQ_p[\overline{G}]\otimes_{\bQ_p}L\cong\bQ_p[\overline{G}]\otimes_{\bQ_p}L^0$$
given by
$$\sigma\otimes m\mapsto\sigma\otimes\sigma^{-1}m,$$
where $L^0$ is the trivial $\overline{G}$-module with $L\cong L^0$ as $\bQ_p-$vector spaces.

When $l\neq p$, we only need to check $\chi(G,M)=-\dim((F^\times/(F^\times)^l)^*\otimes_{\bF_l}\sK)^{\overline{G}}=0$. Discuss it again, we know $((F^\times/(F^\times)^l)^*\otimes_{\bF_l}\sK)^{\overline{G}}\cong((F^\times/\mu\otimes_{\bZ}\bF_l)^*\otimes_{\bF_l}\sK)^{\overline{G}}$. But $F^\times/\mu$ is a $\bZ_p$-module, and $l$ is invertible in $\bZ_p$. Thus, $F^\times/\mu\otimes_{\bZ}\bF_l=0$, $\chi(G,M)=0$.

\begin{remark}: Actually, $dim_{\mathbb{F}_p}{M^{\bar{G}}}=dim_{\mathbb{F}_p}{(M^*)^{\bar{G}}},$
because $\bar{G}$ is a finite cyclic group, let $\sigma$ be the generator of $\bar{G}$, then we have
\begin{eqnarray}
dim_{\mathbb{F}_p}(M^*)^{\bar{G}}&=&dim_{\mathbb{F}_p}(Hom_{\mathbb{F}_p}(M,\mathbb{F}_p))^{\bar{G}}\nonumber\\
&=&dim_{\mathbb{F}_p}Hom_{\mathbb{F}_p[\bar{G}]}(M,\mathbb{F}_p)\nonumber\\
&=&dim_{\mathbb{F}_p}Hom_{\mathbb{F}_p[\bar{G}]}(M/(\sigma-1)M,\mathbb{F}_p)\nonumber\\
&=&dim_{\mathbb{F}_p}Hom_{\mathbb{F}_p}(M/(\sigma-1)M,\mathbb{F}_p)\nonumber\\
&=&dim_{\mathbb{F}_p}M/(\sigma-1)M=dim_{\mathbb{F}_p}M^{\bar{G}},\nonumber
\end{eqnarray}
where M is the $\mathbb{F}_p[\bar{G}]$-module. Thus, by using the conclusion, we know that
\begin{align*}
\dim_{\mathbb{F}_p}{((F^\times/(F^\times)^p)^*\otimes_{\mathbb{F}_p}\sK)^{\bar{G}}}
&=\dim_{\mathbb{F}_p}{(((F^\times/(F^\times)^p)^*\otimes_{\mathbb{F}_p}\sK)^*)^{\bar{G}}}\\
&=\dim_{\mathbb{F}_p}{((F^\times/(F^\times)^p)\otimes_{\mathbb{F}_p}\sK^*)^{\bar{G}}}.
\end{align*}
Thus, we only need to show that
$$\dim_{\mathbb{F}_p}{((F^\times/(F^\times)^p)\otimes_{\mathbb{F}_p}\sK^*)^{\bar{G}}}=[K:\mathbb{Q}_p]\dim_{\mathbb{F}_p}M.$$ Then we can compute it a little easily.
\end{remark}

\section{Global Euler Characteristic Formula}
\begin{theorem}[Global Case]
Let $K/\bQ$ be a finite extension, $S$ be a finite set of places of $\bQ$ including the archimedean places and $\Sigma$ be the set of places of $K$ above $S$. $K^s/K$ is the maximal algebraic extension unramified outside $S$. We write $\fS_S=\gal(K^s/K)$. Assume $M$ is a finite $\fS_S$-module such that if $l\mid |M|$, then $l\in S$. Then we have:
$$\frac{|H^2(\fS_S,M)|\cdot|H^0(\fS_S,M)|}{|H^1(\fS_S,M)|}=\prod_{v\in\Sigma_\infty}\frac{|H^0(G_v,M)|}{||M||_{K_v}},$$
where $|n|_{K_v}=n^{[K_v:\bQ]},$ v is archimedean place, $G_v=\gal(\overline{K}_v/K_v).$
\end{theorem}
In this formula, we know that $M$ is $\fS_S$-module, we also regard $M$ as $G_v$-module.
Before proving this formula, we will give a very powerful theorem which is proved by John Tate. This theorem is the key to prove it. But we do not plan to prove it here. We only narrate it. If you are interested in it ,you can see the reference.
\begin{prop}
Let $S$ be a finite set of places of $\bQ$ including the archimedean places. Let $K$ be a number field and $\Sigma$ be the set of places of $K$ above $S$. Fix a prime $p\in S$, and let $M$ be a discrete finite $\fS_S$-module with $p-$power order. Then we have:
\begin{enumerate}
  \item $$H^r(\fS_S,M)\cong\prod\limits_{v\in\Sigma(\bR)}H^r(G_v,M) \ for \ all \ r\geq 3,$$
where $\Sigma(\bR)$ is the set of real places of $K$.
  \item We have the following long exact sequence:$$0\rightarrow H^0(\fS_S,M)\rightarrow\prod\limits_{v\in\Sigma_0}H^0(G_v,M)\times \prod\limits_{v\in\Sigma_\infty}H_T^0(G_v,M)\rightarrow H^2(\fS_S,M^*(1))^*$$$$\rightarrow H^1(\fS_S,M)\rightarrow\prod\limits_{v\in\Sigma}H^1(G_v,M)\rightarrow H^1(\fS_S,M^*(1))^*$$$$\rightarrow H^2(\fS_S,M)\rightarrow\prod\limits_{v\in\Sigma}H^2(G_v,M)\rightarrow H^0(\fS_S,M^*(1))^*\rightarrow0.$$
\end{enumerate}
\end{prop}

\subsection*{Proof of the Global Case}
Since $M=\bigoplus_lM[l^\infty]$, we may assume $M=M[l^\infty]$. Now, we prove it for $l>2$. Let $$\varphi(M)=\chi(\fS_S,M)-\sum\limits_{v\in\Sigma_\infty}(\len_{\bZ_l}H^0(G_v,M)-[K_v:\bR]\len_{\bZ_l}M),$$ and we need to prove that $\varphi(M)=0$.

By Proposition 2.2(1), we have $H^q(\fS_S,M)=\prod\limits_{v\in\Sigma(\bR)}H^q(G_v,M)$ for all $q\geq3$. However, $\gcd(|G_v|,M)=1$ because of $|G_v|=2$ and $l>3$. Then by Proposition 0.1(2), we have $H^q(G_v,M)=0$ for $q\geq1$, so $H^q(\fS_S,M)=0$ for all $q\geq3$.

When $0\rightarrow L\rightarrow M\rightarrow N\rightarrow0$ is a short exact sequence of $\bZ_l[\fS_S]$-module, we have a long exact sequence:\\
$$0\rightarrow H^0(\fS_S,L)\rightarrow H^0(\fS_S,M)\rightarrow H^0(\fS_S,N)\rightarrow\cdots\rightarrow H^2(\fS_S,N)\rightarrow 0$$\\
of finite $\bZ_l$-module. Hence, $\chi(M)=\chi(L)+\chi(N)$ and $\varphi(M)=\varphi(L)+\varphi(N)$, where $\chi$ and $\varphi$ factor through the Grothendieck group $R_{\bF_l}(\fS_S)$ and have values in $\bZ$.

Then by Proposition 2.2(2), we have $\chi(M)+\chi(M^*(1))=\sum\limits_{v\in\Sigma}\chi_v(G_v,M)$, where
$$\chi_v(G_v,M)=\left\{
                  \begin{array}{ll}
                    \sum\limits_{q=0}^2(-1)^q\len_{\bZ_l}H^q(G_v,M), & v\hbox{ is a non-archimedean place;} \\
                    \sum\limits_{q=0}^2(-1)^q\len_{\bZ_l}H^q_T(G_v,M), & v\hbox{ is an archimedean place.}
                  \end{array}
                \right.
$$Then
$$\sum\limits_{v\in\Sigma}\chi_v(G_v,M)=\sum\limits_{v\in\Sigma_\infty}\chi_v(G_v,M)+\sum\limits_{v\in\Sigma_0}\chi_v(G_v,M),$$
where $\Sigma_0$ is the set of finite places. We know $$\sum\limits_{v\in\Sigma_\infty}\chi_v(G_v,M)=\sum\limits_{v\in\Sigma_\infty}\sum\limits_{q=0}^2(-1)^q\len_{\bZ_l}H^q_T(G_v,M).$$ However, $G_v$ is a cyclic group, by Proposition 0.6, $$\len_{\bZ_l}H^2_T(G_v,M)=\len_{\bZ_l}H^1_T(G_v,M),$$ $$\len_{\bZ_l}H^0_T(G_v,M)=\len_{\bZ_l}H^1_T(G_v,M)=\len_{\bZ_l}H^1(G_v,M).$$ Thus,
$$\sum\limits_{v\in\Sigma_\infty}\chi_v(G_v,M)=\sum\limits_{v\in\Sigma_\infty}\len_{\bZ_l}H^1(G_v,M).$$

On the other hand, since $||M||_{K_v}=1$ for $v\not\in\Sigma$, by the product formula $\prod\limits_{v}||M||_{K_v}=1$ and the local Euler Characteristic formula,
$$\sum\limits_{v\in\Sigma_0}\chi_v(G_v,M)=\log_l\prod_{v\in\Sigma_0}(||M||_{K_v})=-\log_l\prod_{v\in\Sigma_\infty}(||M||_{K_v})=-\sum_{v\in\Sigma_\infty}[K_v:\bR]\len_{\bZ_l}M.$$
Therefore,
\begin{eqnarray*}
\varphi(M)+\varphi(M^*(1))=\chi(M)+\chi(M^*(1))\ \ \ \ \ \ \ \ \ \ \ \ \ \ \ \ \ \ \ \ \ \ \ \ \ \ \ \ \ \ \ \ \ \ \ \ \ \ \ \ \ \ \ \ \ \ \ \ \ \ \ \ \ \ \ \ \ \ \ \ \ \ \ \ \ \ \ \ \ \ \ \ \ \ \ \ \ \ \ \ \ \ \ \ \ \ \ \ \ \ \ \ \ \ \\
-\sum_{v\in\Sigma_\infty}(\len_{\bZ_l}H^0(G_v,M)+\len_{\bZ_l}H^0(G_v,M^*(1))-2[K_v:\bR]\len_{\bZ_l}M)\ \ \ \ \ \ \ \ \ \ \ \ \ \ \ \ \ \ \ \ \ \ \ \ \ \ \ \ \ \ \ \ \ \ \ \ \ \ \ \ \ \ \ \ \ \ \ \ \ \ \ \ \ \ \ \ \ \ \ \ \ \ \ \ \ \ \ \\
=\sum_{v\in\Sigma_\infty}(\len_{\bZ_l}H^1(G_v,M)-\len_{\bZ_l}H^0(G_v,M)-\len_{\bZ_l}H^0(G_v,M^*(1))+[K_v:\bR]\len_{\bZ_l}M)\ \ \ \ \ \ \ \ \ \ \ \ \ \ \ \ \ \ \ \ \ \ \ \ \ \ \ \ \ \ \ \ \ \ \ \ \ \ \ \ \ \ \ \ \ \ \ \ \ \ \ \ \ \ \ \ \ \\
=0. \ \ \ \ \ \ \ \ \ \ \ \ \ \ \ \ \ \ \ \ \ \ \ \ \ \ \ \ \ \ \ \ \ \ \ \ \ \ \ \ \ \ \ \ \ \ \ \ \ \ \ \ \ \ \ \ \ \ \ \ \ \ \ \ \ \ \ \ \ \ \ \ \ \ \ \ \ \ \ \ \ \ \ \ \ \ \ \ \ \ \ \ \ \ \ \ \ \ \ \ \ \ \ \ \ \ \ \ \ \ \ \ \ \ \ \ \ \ \ \ \ \ \ \ \ \ \ \ \ \ \ \ \ \ \ \ \ \ \ \ \ \ \ \ \ \ \ \ \ \ \ \ \ \ \ \ \ \ \ \ \ \ \ \ \ \ \ \ \ \ \ \ \ \ \ \ \ \ \ \
\end{eqnarray*}
The last equation follows the following proposition:
\begin{prop}
(in [2, theorem2.3(c)]) Let $K_v=\bR$ or $\bC$ and let $G_v=\gal(\overline{K}_v/K_v)$, $|n|_{K_v}=n^{[K_v:\bR]}$. For any finite $G_v$-module $M$, Then we have
$$\frac{|H^0(G_v,M))|\cdot|H^0(G_v,M^*(1))|}{|H^1(G_v,M)|}=||M||_{K_v}.$$
\end{prop}
So we only need to show that $\varphi(M)=\varphi(M^*(1))$.

We take a finite Galois extension $F/K$ such that $\fS_S'=\gal(K^s/F)$ acts on $M$ and $\mu_l$ trivially. We write $\overline{G}=\gal(F/K)$. By the same argument of local case, we may assume $\overline{G}$ is cyclic of degree prime to $l$ and $M$ is a $\bF_l[\overline{G}]$-module. We know $H^q(\fS_S',M)=0$ for $q\geq3$, since $\gcd(|\overline{G}|,l)=1$, $H^q(\overline{G},M)=0$ for all $q>0$, and by the inflation and restriction sequence again, we get $H^q(\fS_S,M)\cong H^0(\overline{G},H^q(\fS_S',M))$.

Since we consider $\varphi$ as a homomorphism from the Grothendieck group $R_{\bF_l}[\overline{G}]$ to $\bQ$. let $$\chi':R_{\bF_l}(\overline{G})\otimes\bQ\rightarrow R_{\bF_l}(\overline{G})\otimes\bQ,$$$$[M]\mapsto\sum\limits_{i=0}^2(-1)^i[H^i(\fS_S',M)]$$ and $$\theta:R_{\bF_l}(\overline{G})\otimes\bQ\rightarrow\bQ,$$$$[M]\mapsto\dim_{\bF_l}M^{\overline{G}},$$ then $$\chi=\theta\circ\chi'.$$
We know that $[H^0(\fS_S',\mu_l)]=[\mu_l]$, and
\begin{description}
  \item[(i)] $H^1(\fS_S',\mu_l)=[O_{F,S}^\times/l]+[Cl_S(F)[l]]$,
  \item[(ii)] $[H^2(\fS_S',\mu_l)]=[Cl_S(F)/l]-[\bF_l]+[\bigoplus\limits_{\mathfrak{p}\in S\backslash S_\infty(F)}\bF_l]+[\bigoplus\limits_{\mathfrak{p}\in S_\infty(F)}H^0_T(G_\mathfrak{p},\bF_l)]$,
\end{description}
where $O^\times_{F,S}$ is the group of $S-$unit and $Cl_S(F)[l]$ is the $l-$torsion part of $S-$ideal class group $Cl_S(F)$. See more details in [8, (8.7.4)].

We obtain that $[O_{F,S}^\times/l]=[\bigoplus\limits_{\mathfrak{p}\in S(F)}\bF_l]+[\mu_l]-[\bF_l]$ and $[Cl_S(F)[l]]=[Cl_S(F)/l]$, then we have $$\chi'([\mu_l])=[\bigoplus\limits_{\mathfrak{p}\in S_\infty(F)}H^0_T(G_\mathfrak{p},\bF_l)]-[\bigoplus\limits_{\mathfrak{p}\in S_\infty(F)}\bF_l].$$ For a finite $\bF_l[\overline{G}]$-module $M$, we have
\begin{description}
  \item[(iii)] $\chi'(M^*(1))=[M^*]\cdot\chi'([\mu_l])$, see [2, Lemma 5.4],
  \item[(iv)] $[M]\cdot[\bF_l[\overline{G}]]=\dim_{\bF_l}M\cdot[\bF_l[\overline{G}]]$
\end{description}
Thus
\begin{eqnarray*}
\chi'([M^*(1)])&=&[M^*\otimes(\bigoplus\limits_{\mathfrak{p}\in S_\infty(F)}(H^0_T(G_\mathfrak{p},\bF_l)-\bF_l))]\\
&=&[M^*\otimes(\bigoplus\limits_{v\in\Sigma_\infty}\bigoplus\limits_{\mathfrak{p}|v}(H^0_T(G_\mathfrak{p},\bF_l)-\bF_l))].
\end{eqnarray*}
We obtain that if $v=\bC$, $\mathfrak{p}=\bC$, then $$\bigoplus\limits_{\mathfrak{p}|v}\bF_l\cong\bF_l[\overline{G}]$$ as $\overline{G}$-modules and $$H^0_T(G_\mathfrak{p},\bF_l)=0.$$

If $v=\bR$, $\mathfrak{p}=\bR$, then $$\bigoplus\limits_{\mathfrak{p}|v}\bF_l\cong\bF_l[\overline{G}]$$ as $\overline{G}$-modules and
$$\bigoplus\limits_{\mathfrak{p}|v}H^0_T(G_\mathfrak{p},\bF_l)\cong H^0_T(G_\mathfrak{p},\bF_l)\otimes_{\bF_l}\bF_l[\overline{G}].$$
If $v=\bR$, $\mathfrak{p}=\bC$, then $H^0_T(G_\mathfrak{p},\bF_l)=0$.

We can compute similarly for $\chi'(M)$. Then by (iv), for $\varphi(M)=\varphi(M^*(1))$, we only need to prove that
\begin{eqnarray*}
-\theta([M^*\otimes(\bigoplus\limits_{v\in\Sigma_\infty(\bR)}\bigoplus\limits_{\mathfrak{p}|v,\mathfrak{p}\hbox{ is complex}}\bF_l)])
-\sum\limits_{v\in\Sigma_\infty}(\dim_{\bF_l}H^0(G_v,M^*(1))-[K_v:\bR]\dim_{\bF_l}M^*(1))\\
=-\theta([M(-1)\otimes(\bigoplus\limits_{v\in\Sigma_\infty(\bR)}\bigoplus\limits_{\mathfrak{p}|v,\mathfrak{p}\hbox{ is complex}}\bF_l)])
-\sum\limits_{v\in\Sigma_\infty}(\dim_{\bF_l}H^0(G_v,M)-[K_v:\bR]\dim_{\bF_l}M),
\end{eqnarray*}
where $M(-1)=(M^*(1))^*$.

In other words, we only need to prove
\begin{eqnarray*}
&\sum\limits_{v\in\Sigma(\bR)}(\theta([M^*\otimes(\bigoplus_{\mathfrak{p}|v,\mathfrak{p}\mbox{ is complex}}\bF_l)])+\dim_{\bF_l}M^*(1)^{G_v})\\
=&\sum\limits_{v\in\Sigma(\bR)}(\theta([M(-1)\otimes(\bigoplus_{\mathfrak{p}|v,\mathfrak{p}\mbox{ is complex}}\bF_l)])+\dim_{\bF_l}M^{G_v}).
\end{eqnarray*}
We will prove the equation for each $v$.

However, $$\theta([M^*\otimes(\bigoplus\limits_{\mathfrak{p}|v,\mathfrak{p}\mbox{ is complex}}\bF_l)])=\dim_{\bF_l}(M^*\otimes(\bigoplus\limits_{\mathfrak{p}|v,\mathfrak{p}\mbox{ is complex}}\bF_l))^{\overline{G}}.$$

We know that $G_v$ is the subgroup of $\overline{G}$, define $\overline{G}'=\overline{G}/G_v$, then we have
\begin{eqnarray*}
\theta([M^*\otimes(\bigoplus\limits_{\mathfrak{p}|v,\mathfrak{p}\mbox{ is complex}}\bF_l)])
&=&\dim_{\bF_l}((M^*\otimes(\bigoplus\limits_{\mathfrak{p}|v,\mathfrak{p}\mbox{ is complex}}\bF_l))^{G_v})^{\overline{G}'}\\
&=&\dim_{\bF_l}((M^*)^{G_v}\otimes(\bigoplus\limits_{\mathfrak{p}|v,\mathfrak{p}\mbox{ is complex}}\bF_l))^{\overline{G}'}\\
&=&\dim_{\bF_l}((M^*)^{G_v}\otimes\bF_l[\overline{G}'])^{\overline{G}'}\\
&=&\dim_{\bF_l}(M^*)^{G_v}\\
&=&\dim_{\bF_l}M^{G_v}.
\end{eqnarray*}
Where the last two equations are separately given by (iv) and remark(1.5).
By the same discussion, we get
$$\theta([M(-1)\otimes(\bigoplus\limits_{\mathfrak{p}|v,\mathfrak{p}\mbox{ is complex}}\bF_l)])=\dim_{\bF_l}M^*(1)^{G_v}.$$
Hence, we complete the proof for $l>2$.

When l=2, we have $M^*(1)=M^*$ and we need to modify the proof of additions of $\varphi$ which could be seen in [1, 4.4.5].
\begin{remark}
They directly prove $\chi(M)=\chi(M^*(1))$ and $\dim_{\bF_l}M^{G_v}=\dim_{\bF_l}(M^*(1))^{G_v}$ in [1,Ch4.4.5] and [2,Ch1.5]. Actually, we cannot get $\chi(M)=\chi(M^*(1))$ directly from $\chi'$. And also we cannot get $\dim_{\bF_l}M^{G_v}=\dim_{\bF_l}(M^*(1))^{G_v}$. A counterexample is that let $K=\mathbb{Q}(\sqrt{3})$ and $L=K(\sqrt{-1})=\mathbb{Q}(\sqrt{3},\sqrt{-1})\supset \mu_3$ be a cyclic extension of $K$ of degree $2$
which primes to $3$. The prime in $K$ above $3$ is unramified in $L$. Let $M=\mathbb{F}_3$, then $M^*(1)=\mu_3$. But
$\dim_{\bF_l}M^{G_v}=3\neq 1=\dim_{\bF_l}(M^*(1))^{G_v}.$
\end{remark}

\end{document}